\documentclass[12pt,leqno]{article}

\usepackage{amsmath,amssymb,amsthm}

\makeatletter

\newtheorem{thm}{Theorem}[section]
\newtheorem{prp}[thm]{Proposition}
\newtheorem{lmm}[thm]{Lemma}
\newtheorem{crl}[thm]{Corollary}

\theoremstyle{definition}

\theoremstyle{remark}
\newtheorem{rem}[thm]{Remark}

\def\({{\rm (}}
\def\){{\rm )}}

\let\Mathrm\operator@font
\let\Cal\mathcal
\newcommand{\fm}{\ensuremath{\mathfrak m}}

\def\standop#1{\mathop{\Mathrm #1}\nolimits}
\def\difstop#1#2{\expandafter\def\csname #1\endcsname{\standop{#2}}}
\def\defstop#1{\difstop{#1}{#1}}

\defstop{Hom}
\defstop{Gal}
\defstop{Spec}

\defstop{Ker}
\defstop{Min}
\difstop{Image}{Im}
\difstop{tdeg}{trans.deg}
\def\red{_{\Mathrm{red}}}
\defstop{length}
\defstop{supp}
\difstop{height}{ht}
\def\uHom{\mathop{\mbox{\underline{$\Mathrm Hom$}}}\nolimits}
\def\uann{\mathop{\mbox{\underline{$\Mathrm ann$}}}\nolimits}

\def\specialarrow#1{\setbox\z@=\hbox{$\m@th
 \mathop{\vphantom{\rightarrow}}\limits^{\hspace{.5ex}{#1}\hspace
{.8ex}}$}\mathrel{\ifdim\wd\z@<1.2em\dimen\tw@
1.2em\else\dimen\tw@\wd\z@\fi\copy\z@\kern-\wd\z@\hbox to\dimen\tw@
{\rightarrowfill}}}

\def\sdarrow#1{\downarrow\hbox to 0pt{\scriptsize$#1$\hss}}
\def\suarrow#1{\uparrow\hbox to 0pt{\scriptsize$#1$\hss}}

\def\section{\@startsection{section}{1}{\z@ }%
{-3.5ex plus -1ex minus -.2ex}{2.3ex plus .2ex}{\bf }}

\long\def\refname{\par\kern -3ex
\begin{center}\rm R\sc{eferences}\end{center}\par\kern 
-2ex}

\def\@seccntformat#1{\csname the#1\endcsname.\quad}

\def\@@@sect#1#2#3#4#5#6[#7]#8{%
   \ifnum #2>\c@secnumdepth 
      \def \@svsec {}\else \refstepcounter {#1}%
      \def\@svsec{}
   \fi 
   \@tempskipa #5\relax 
   \ifdim \@tempskipa >\z@ 
     \begingroup #6\relax \@hangfrom {\hskip #3\relax 
     \@svsec}{\interlinepenalty \@M #8\par }\endgroup 
     \csname #1mark\endcsname {#7}
   \else 
   \def \@svsechd {#6\hskip #3\@svsec #8\csname #1mark\endcsname {#7}}
   \fi \@xsect {#5}}

\def\@@@startsection#1#2#3#4#5#6{%
 \if@noskipsec \leavevmode \fi \par \@tempskipa #4\relax \@afterindenttrue 
 \ifdim \@tempskipa <\z@ \@tempskipa -\@tempskipa \@afterindentfalse 
 \fi \if@nobreak \everypar {}\else \addpenalty {\@secpenalty }\addvspace 
  {\@tempskipa }\fi \@ifstar {\@ssect {#3}{#4}{#5}{#6}}{\@dblarg 
  {\@@@sect {#1}{#2}{#3}{#4}{#5}{#6}}}}

\def\theparagraph{\thesection.\arabic{paragraph}}
\def\aparagraph{\@@@startsection{paragraph}{2}{\z@ }%
              {1.75ex plus .2ex minus .15ex}{-1em}{\bf(\theparagraph) } }
\def\paragraph{\@@@startsection{paragraph}{2}{\z@ }%
              {1.75ex plus .2ex minus .15ex}{-1em}{}{\bf(\theparagraph)} }

\let\c@Theorem\c@paragraph

\title{``Geometric quotients are algebraic schemes'' based on Fogarty's 
idea}

\author{M{\sc itsuyasu} H{\sc ashimoto}}
\date{\normalsize
Graduate School of Mathematics, Nagoya University\\
Chikusa-ku,  Nagoya 464--8602 JAPAN\\
\small \tt hasimoto@math.nagoya-u.ac.jp}

\begin{document}

\maketitle

\begin{abstract}
Let $S$ be a Noetherian scheme, $\varphi:X\rightarrow Y$ a surjective
$S$-morphism of $S$-schemes, with $X$ of finite type over $S$.
We discuss what makes $Y$ of finite type.

First, we prove that if $S$ is excellent, $Y$ is reduced, and
$\varphi$ is universally open, then $Y$ is of finite type.
We apply this to understand Fogarty's theorem in ``Geometric quotients are 
algebraic schemes, Adv.\ Math.\ {\bf 48} (1983), 166--171''
for the special case that the group scheme $G$ is flat over the
Noetherian base scheme $S$ and that the quotient map is universally
submersive.
Namely, we prove that if $G$ is a flat $S$-group scheme of finite type
acting on $X$ and $\varphi
$ is its universal strict orbit space, then $Y$ is of finite type
($S$ need not be excellent.
Geometric fibers of $G$ can be disconnected and non-reduced).

Utilizing the technique used there, we also prove that $Y$ is of finite type
if $\varphi$ is flat.
The same is true if $S$ is excellent, $\varphi$ is proper, 
and $Y$ is Noetherian.
\end{abstract}

\section{Introduction}

In \cite{Fogarty}, Fogarty proved the following theorem.

\begin{thm}[Fogarty]
Let $S$ be an excellent scheme, $G$ an $S$-group scheme of finite type
with connected geometric fibers.
Let $X$ be a $G$-scheme of finite type over $S$.
If $(Y,\varphi)$ is a strict orbit space for the action of $G$ on $X$, then
$Y$ is of finite type over $S$.
Moreover, if $\Cal F$ is a coherent $(G,\Cal O_X)$-module
\(coherent $G$-linearlized $\Cal O_X$-module\),
then $(\varphi_*\Cal F)^G$ is a coherent $\Cal O_Y$-module.
\end{thm}

The main purpose of this paper is to try to 
understand this very important theorem
in invariant theory.
The author has not understood his proof yet.
On the other hand, his idea is transparent, and
if we assume that $G$ is $S$-flat and 
$\varphi$ is universally submersive, then one can keep track of his proof
without much difficulty, and even remove some other assumptions.
Namely, we prove

\begin{thm}\label{geometric-quotient.thm}
Let $S$ be a Noetherian scheme, $G$ a flat $S$-group scheme of finite type.
Let $X$ be a $G$-scheme of finite type over $S$.
If $(Y,\varphi)$ is a universal 
strict orbit space for the action of $G$ on $X$, then
$Y$ is of finite type over $S$.
Moreover, if $\Cal F$ is a coherent $(G,\Cal O_X)$-module,
then $(\varphi_*\Cal F)^G$ is a coherent $\Cal O_Y$-module.
\end{thm}

Flatness assumption is important in our proof, since
we use the assumption that the image and the kernel of a $G$-equivariant
$\Cal O_X$-module map between $(G,\Cal O_X)$-modules 
are again $(G,\Cal O_X)$-modules.
Another merit in using flatness is the universal openness of 
universal orbit spaces \cite[p.6]{GIT}.

We do not require that $S$ is excellent.
We do not require that $G$ has connected geometric fibers.
The theorem includes the case that $X=G$.
So we do not claim that $X\red$ is $G$-stable, or that irreducible components
of $X$ are $G$-stable.

The proof for the most essential case, the case 
that $Y$ is reduced and $S$ is excellent 
is purely ring-theoretic, see Theorem~\ref{main.thm}.
The proof heavily depends on the idea of Fogarty \cite{Fogarty} and
Onoda \cite{Onoda}.
To remove excellence assumption of the base $S$, we utilize 
Onoda's result \cite[(2.20)]{Onoda}.

Utilizing the technique used above, 
we also prove some finite generation results which do not have direct
connection to group actions.
Let $S$ be a Noetherian scheme, and 
$\varphi:X\rightarrow Y$ a surjective morphism of $S$-schemes, with $X$
of finite type.
If $\varphi$ is flat, then $Y$ is of finite type (Corollary~\ref{flat.thm}).
If $S$ is excellent, 
$\varphi$ is proper, and $Y$ is Noetherian, then $Y$ is of finite type
(Theorem~\ref{proper.thm}).

Onoda \cite{Onoda} proved that if $S$ is Nagata and all
normal local rings that are essentially of finite type over $S$ are
analytically irreducible (e.g., $S$ is excellent), $Y$ is Noetherian 
normal, and the generic point of any irreducible component of $X$ is
mapped to the generic point of an irreducible component of $Y$, 
then $Y$ is of finite type.
Fogarty \cite{Fogarty2} proved the same result independently later.
Our argument more or less follows theirs,
but none of our new assertions here does not cover their theorem.

The author is grateful to Professor Shigeru Mukai, Professor Tetsushi
Ogoma and Professor Nobuharu Onoda for valuable advice.
Special thanks are also due to the referee for valuable comments.

\section{Main theorem --- the reduced case}

Throughout this section, let $S$ be a Noetherian scheme, 
$\varphi:X\rightarrow Y$ a surjective $S$-morphism of $S$-schemes,
with $X$ of finite type over $S$.
The following proposition is in \cite{Fogarty}.
We state and give a proof for readers' convenience.

\begin{prp}\label{prop.thm}
Assume that $S=\Spec R$ is affine, and 
$Y=\Spec B$ is affine with $B$ an integral domain.
Then 
there is a finitely generated $R$-subalgebra $A$ of $B$ such that
the induced morphism $\eta:Y=\Spec B\rightarrow \Spec A=:Z$ is 
birational and geometrically injective.
\end{prp}

\begin{proof} Let $(U_i)$ be a finite affine open covering of $X$.
Then replacing $X$ by $\coprod_i U_i$, we may assume that
$X=\Spec C$ is affine.

Let $J$ be the kernel of the canonical map $C\otimes_R C
\rightarrow C\otimes_B C$.
Then there exist some $b_1,\ldots,b_r\in B$ such that 
$J$ is generated by
\[
b_1\otimes 1-1\otimes b_1,\ldots,
b_r\otimes 1-1\otimes b_r.
\]
Set $A_0=R[b_1,\ldots,b_r]$.
Let $y, y'\in Y(\xi)$ be distinct geometric points of $Y$,
where $\xi$ is an algebraically closed field.
Since $\varphi$ is surjective of finite type, there exist $x,x'\in X(\xi)$
such that $\varphi(x)=y$ and $\varphi(x')=y'$.
Since $(x,x')\in (X\times X)(\xi)\setminus(X\times_Y X)(\xi)$, 
there exists some
$i$ such that $b_i(x)\neq b_i(x')$.
This shows that the image of $y$ and $y'$ in $(\Spec A_0)(\xi)$ are 
different.
So $Y\rightarrow \Spec A_0$ is geometrically injective.
It is clear that there exist some $a_1,\ldots,a_t\in B$ such that
$A=A_0[a_1,\ldots,a_t]$ is birational to $B$.
Then $A$ is the desired subalgebra, since $Y\rightarrow \Spec A$ is
still geometrically injective.
\end{proof}

The following is \cite[Lemma~3]{Fogarty}.

\begin{lmm}\label{lemma.thm}
Let $\psi: U\rightarrow Z$ be an affine birational morphism between 
integral schemes.
If $Z$ is Noetherian 
normal and $\psi(U)$ is an open subset of $Z$, then $\psi$ is
an open immersion.
\end{lmm}

We omit the proof.
\iffalse
\proof Clearly, $\psi$ factors through the open subscheme 
$\psi(U)$ of $Z$.
So we may assume that $\psi$ is surjective.
As the question is local, we may assume that $Z=\Spec A$ is affine.
As $\psi$ is affine, $U=\Spec B$ is affine.
By assumption, we have $A\subset B\subset Q(A)$, where $Q(A)$ is the 
field of fractions of $A$.

Let $P$ be a height one prime ideal of $A$.
Then $A_P\subset B_P\subset Q(A)$, and there is a prime ideal of $B_P$
lying over $PA_P$.
This shows $A_P=B_P$, otherwise we must have $B_P=Q(A)$.
So $B\subset B_P\subset A_P$ for any $P$, and it follows that 
$B\subset \bigcap_P A_P=A$.
\qed
\fi

The following is based on the ideas of Fogarty \cite{Fogarty}
and Onoda \cite{Onoda}.

\begin{thm}\label{main.thm}
Let $S$ be a universally catenary Nagata scheme
\(i.e., for any affine open subset $U=\Spec R$ of $S$, $R$ is
universally catenary and Nagata\), and $\varphi
:X\rightarrow Y$ a surjective 
universally open $S$-morphism of $S$-schemes.
If $X$ is of finite type over $S$ and $Y$ is reduced, 
then $Y$ is of finite type over $S$.
\end{thm}

\begin{proof} Clearly, $Y$ is quasi-compact.
So the question is local on $S$ and $Y$, and so 
we may assume that $S=\Spec R$ and $Y=\Spec B$ are
affine.
Let $(U_i)$ be a finite affine open covering of $X$.
Then replacing $X$ by $\coprod_i U_i$, we may assume that $X=\Spec C$ is
also affine.

Since $X$ has only finitely many irreducible components and $\varphi$ is
surjective, $B$ has only finitely many minimal primes.
Since $B\rightarrow \prod_{P\in\Min(B)}B/P$ is injective and finite,
it suffices to show that $B/P$ is of finite type for $P\in\Min(B)$.
Replacing $B$ by $B/P$, we may assume that $B$ is a domain.

Take $A\hookrightarrow B$ as in Proposition~\ref{prop.thm}
so that $\eta:Y=\Spec B\rightarrow \Spec A=Z$ is geometrically injective
and birational,
and $Z$ is of finite type over $R$.
Now 
let $A'$ be the normalization of $A$, and $B'=B[A']$.
Since $R$ is Nagata, the associated morphism $\alpha:Z'=\Spec A'\rightarrow
\Spec A=Z$ is finite.
Let $Y'=\Spec B'$ and $X'=X\times_Y Y'$.
\[
\begin{array}{ccccc}
X' & \specialarrow{\varphi'} & Y' & \specialarrow{\eta'} & Z'\\
\sdarrow{\gamma} &           & \sdarrow{\beta} &    & \sdarrow{\alpha}\\
X  & \specialarrow{\varphi}  & Y  & \specialarrow{\eta} & Z
\end{array}
\]

Note that $Y'$ is a closed subscheme of $Y\times_Z Z'$.
In particular, $\beta$ is finite, since $\alpha$ and the closed immersion are.
Similarly, $\eta'$ is geometrically injective.
Clearly, $\gamma$ is finite and $\varphi'$ is universally open.

Let $x'\in X'$.
Set $y'=\varphi'(x')$, $z'=\eta'(y')$, 
$F=\Cal O_{X',x'}$, $E=\Cal O_{Y',y'}$, and $D=\Cal O_{Z',z'}$.
Then for any minimal prime $Q$ of $F$, we have
$Q\cap E=0$, since $\varphi'$ is open.

Since $Z'$ is universally catenary, 
there exists some $n\geq 0$ such that
\begin{multline*}
\dim E(t_1,\ldots,t_n)-\dim D
=\\
\tdeg_{R(Z')}R(Y')(t_1,\ldots,t_n)-
\tdeg_{\kappa(z')}\kappa(y')(t_1,\ldots,t_n)
\end{multline*}
by Onoda's dimension formula \cite[(1.11)]{Onoda}, 
where $t_1,\ldots,t_n$ are variables,
and for a local ring $(O,\frak m)$, $O(t_1,\ldots,t_n)$ denotes the
local ring $O[t_1,\ldots,t_n]_{\frak m[t_1,\ldots,t_n]}$.
Since $\eta'$ is universally injective, $\kappa(y')$ is a purely 
inseparable algebraic extension of $\kappa(z')$
by \cite[(3.5.8)]{EGA-I}.
Thus the right hand side is zero, since $R(Z')=R(Y')$.

Let $P$ be a minimal prime of $\fm_{y'}F$ such that 
$\dim F\otimes_E \kappa(y')=\dim F/P$.
Since $\varphi'\times 1:X'\times \Bbb A^n \rightarrow Y'\times \Bbb A^n$ 
is an open map, we have
\[
\height P=\height P[t_1,\ldots,t_n]\geq \dim E(t_1,\ldots,t_n)=\dim D.
\]
So we have
\[
\dim F\geq \dim D+\dim (F\otimes_E\kappa(y'))
=\dim D+\dim(F\otimes_D\kappa(z'))
\]
by the geometric injectivity of $\eta'$.
Hence we have $\dim F=\dim D+\dim(F\otimes_D\kappa(z'))$
by \cite[(15.1)]{CRT}.

For $r\geq 0$, define $X'(r)$ to be $(\bigcup X_r'
)\setminus(\bigcup X_{>r}')$,
where $X_r'$ (resp.\ $X_{>r}'$) runs through all irreducible components
(with reduced structures) 
of $X'$ such that $\tdeg_{R(Y')}R(X_r')=r$ (resp.\ 
$\tdeg_{R(Y')}R(X_{>r}')>r$).
Let $x'\in X'(r)$.
By the dimension formula \cite[(14.C)]{CA}, 
we have that
\begin{multline*}
\dim \Cal O_{X'(r),x'}=
\dim \Cal O_{X',x'}=
\dim \Cal O_{Z',z'}+\dim \Cal O_{X',x'}\otimes_{\Cal O_{Z',z'}}\kappa(z')\\
\geq
\dim \Cal O_{Z',z'}+\dim \Cal O_{X'(r),x'}\otimes_{\Cal O_{Z',z'}}\kappa(z'),
\end{multline*}
where $z'=(\eta'\varphi')(x')$.
So we have 
\[
\dim \Cal O_{X'(r),x'}=
\dim \Cal O_{Z',z'}+\dim \Cal O_{X'(r),x'}\otimes_{\Cal O_{Z',z'}}\kappa(z').
\]
Since all local rings of $X'(r)$ are equidimensional by the dimension 
formula, we have that $\eta'\varphi'|_{X'(r)}$ is equidimensional for
any $r$ \cite[(13.3.6)]{EGA-IV}.
Since $Z'$ is normal, $\eta'\varphi'|_{X'(r)}$ is universally open 
for any $r$, by Chevalley's criterion \cite[(14.4.4)]{EGA-IV}.
Since $X'=\bigcup_{r\geq 0} X'(r)$, we have that $\eta'\varphi'$ is
universally open.

Since $\varphi'$ is surjective, $\eta'(Y')=(\eta'\varphi')(X')$ is open
in $Z'$.
Since $\eta'$ is affine birational and $Z'$ is Noetherian normal,
$\eta'$ is an open immersion by Lemma~\ref{lemma.thm}.
So $Y'$ is of finite type.
Since $B'$ is finitely generated and 
$B\rightarrow B'$ is finite and injective, $B$ is finitely generated.
\end{proof}

\begin{rem}
We can also prove the following.
Let $S$ be a Nagata scheme, and $\varphi:X\rightarrow Y$ a surjective 
universally open equidimensional $S$-morphism of $S$-schemes.
If $X$ is of finite type over $S$ and $Y$ is reduced, 
then $Y$ is of finite type over $S$.
Indeed, 
since $\varphi$ is universally open and equidimensional, $\varphi'$ is so.
Since $\varphi'$ is equidimensional and $\eta'$ is geometrically injective
and birational, 
$\eta'\varphi'$ is also equidimensional, and the same argument works.
\end{rem}

\begin{crl}\label{reduced.thm}
Let $S$ be an excellent scheme, $G$ a flat $S$-group scheme of finite type,
and $X$ a $G$-action of finite type over $S$.
If $\varphi:X\rightarrow Y$ is a universal 
strict orbit space and $Y$ is reduced,
then $Y$ is of finite type.
\end{crl}

\begin{proof} $\varphi$ is universally open and 
surjective, see \cite[p.6]{GIT}.
So we are done.
\end{proof}

\begin{crl}\label{flat.thm}
Let $S$ be a Noetherian scheme, and $\varphi:X\rightarrow Y$ 
a faithfully flat $S$-morphism of $S$-schemes.
If $X$ is of finite type over $S$, then $Y$ is of finite type over $S$.
\end{crl}

\begin{proof} We may assume that $S=\Spec R$, $Y=\Spec B$, and $X=\Spec C$ are
all affine.
Since $B$ is a pure subring of $C$ and $C$ is Noetherian, we have that
$B$ is also Noetherian.
Since the nilradical of $B$ is a finitely generated ideal, it suffices to
show that $B\red$ is of finite type.
Replacing $B$ by $B\red$ and $C$ by $C\otimes_B B\red$, we may assume
that $B$ is reduced.
So it suffices to show that $\prod_{P\in \Min(B)} B/P$ is of finite type.
Replacing $B$ by $B/P$, we may assume that $B$ is a domain.
As $C$ is faithfully flat over the domain $B$, there exists some
prime ideal $Q$ of $C$ such that $Q\cap B=0$.
By \cite[(2.11) and (2.20)]{Onoda}, we may assume that $R$ is local.
By \cite[(2.7.1)]{EGA-IV}, we may assume that $R$ is a complete local ring.
Replacing $B$ again if necessary, we may still assume that $B$ is a domain.
Since $R$ is excellent and $B$ is a domain, 
$B$ is of finite type by the theorem.
\end{proof}

\section{The general case}

In this section, we prove Theorem~\ref{geometric-quotient.thm}.
First we prove the following

\begin{lmm}\label{reduction-lemma.thm}
Let $S$ be a Noetherian scheme, and $G$ a flat $S$-group scheme of finite type.
Let $X$ be a $G$-scheme of finite type over $S$.
Assume that for any closed $G$-subscheme $X_1$ of $X$ and its
universal
strict orbit space $\psi: X_1\rightarrow Y_1$, we have that $Y_1$ is
of finite type, provided $Y_1$ is reduced \(e.g., $S$ is excellent, see 
{\rm Corollary~\ref{reduced.thm}}\).
If $(Y,\varphi)$ is a universal
strict orbit space for the action of $G$ on $X$, then
$Y$ is of finite type over $S$.
Moreover, if $\Cal F$ is a coherent $(G,O_X)$-module,
then $(\varphi_*\Cal F)^G$ is a Noetherian $\Cal O_Y$-module.
\end{lmm}

\begin{proof} We may assume that $S=\Spec R$ and $Y=\Spec B$ are affine.

First, we prove the last assertion using the Noetherian induction 
on the coherent 
ideal sheaf $\uann\Cal F:=\Ker(\Cal O_X\rightarrow \uHom_{\Cal O_X}
(\Cal F,\Cal F))$.
We may assume that $\Cal F\neq 0$.
We also use the induction on 
$\nu(\Cal F):=\sum_V \length_{\Cal O_{X,v}}\Cal F_v$,
where $V$ runs through the irreducible components of $\supp\Cal F
=V(\uann\Cal F)$, and $v$ is the generic point of $V$.

Let $\Cal G$ be the maximal coherent $(G,\Cal O_X)$-submodule 
of $\Cal F$ such that 
$(\varphi_*\Cal G)^G$ is a Noetherian $\Cal O_Y$-module.
As $(\varphi_*?)^G$ is left exact, we may assume that $\Cal G=0$, 
replacing $\Cal F$ by $\Cal F/\Cal G$.
In particular, any nonzero coherent $(G,\Cal O_X)$-subsheaf of $\Cal F$
has the same annihilator as that of $\Cal F$.

If $(\varphi_*\Cal F)^G=0$, then there is nothing to be proved.
So we consider the case that
$H^0(X,\Cal F)^G=(\varphi_*\Cal F)^G\neq 0$.
Take $a\in H^0(X,\Cal F)^G\setminus\{0\}=\Hom_{G,\Cal O_X}(\Cal O_X,\Cal F)
\setminus \{0\}$.
Then
$a\Cal O_X$ is a nonzero coherent $(G,\Cal O_X)$-subsheaf of $\Cal F$.
So $\nu(a\Cal O_X)\neq 0$.
If $\nu(\Cal F)>\nu(a\Cal O_X)$, then $\Cal G\supset a\Cal O_X\neq 0$
by induction assumption, and this is a contradiction.
So $\nu(\Cal F)=\nu(a\Cal O_X)$.
So $\nu(\Cal F/a\Cal O_X)=0$, which shows that 
$\supp(\Cal F/a\Cal O_X)\subsetneq \supp\Cal F$.
By induction assumption, $\varphi_*(\Cal F/a\Cal O_X)^G$ is Noetherian.
So $\varphi_*(a\Cal O_X)^G$ is not Noetherian, since we assume that
$\Cal F\neq 0=\Cal G$.
Let $\Cal J$ be the kernel of $a:\Cal O_X\rightarrow\Cal F$, and
$Z$ be the $G$-stable closed subscheme of $X$ defined by the
coherent $G$-ideal $\Cal J$ of $\Cal O_X$.
We may assume that $\Cal F=\Cal O_Z$, from the beginning.
If $H^0(X,\Cal O_Z)^G$ is not reduced, then there exists some
$b\in H^0(X,\Cal O_Z)^G\setminus\{0\}$ such that $b^2=0$.
Then $0\neq b\Cal O_Z\subset \Cal O_Z$, and the annihilator of
$b\Cal O_Z$ is strictly larger than that of $\Cal O_Z$.
By induction assumption, this is a contradiction.
So $B_1=H^0(X,\Cal O_Z)^G$ must be a reduced ring.

Set $Y_1:=\Spec B_1$, $B_0$ to be the image of the canonical map
$B\rightarrow B_1$, 
and $Y_0=\Spec B_0$ the scheme theoretic image of 
$Y_1\rightarrow Y$.
Let $\varphi_1:Z\rightarrow Y_1$ and $\eta':Y_1\rightarrow Y_0$ 
be the canonical maps.
Note that $Y_0=\varphi(Z)=\Image(\eta'\varphi_1)$ set-theoretically, 
since $Z$ is $G$-stable closed, and $\varphi$ is an orbit space.
We have $\varphi^{-1}(Y_0)=Z$ set theoretically, so 
the inclusion $Z\hookrightarrow\varphi^{-1}(Y_0)$ is a universal 
homeomorphism.
Hence $\eta'\varphi_1$ is surjective and universally open.
It follows that $(Y_0,\eta'\varphi_1)$ is a universal strict orbit space.
Since $Y_0$ is reduced, $Y_0$ is of finite type by assumption.

By \cite[Proposition~1]{Fogarty}, $(Y_1,\varphi_1)$ is a universal 
geometric 
quotient, and $\eta'$ is a universal homeomorphism.
In particular, $Y_1$ is of finite type by assumption.
So $\eta':Y_1\rightarrow Y_0$ is a universal homeomorphism
of finite type between Noetherian schemes.
So $\eta'$ is finite.
Hence $B_1$ is a $B_0$-finite module.
Since $B_0$ is of finite type over $R$, $B_1$ is a Noetherian $B_0$-module.
Hence $(\varphi_*\Cal O_Z)^G$ is a Noetherian $\Cal O_Y$-module,
as desired.

Next, we prove that $Y=\Spec B$ is of finite type.
We use Noetherian induction, and we may assume that for a $G$-stable
closed subscheme $X_1\subsetneq X$ and its universal strict orbit space
$\varphi_1:X_1\rightarrow Y_1$, $Y_1$ is of finite type.
If $Y$ is reduced, then there is nothing to be proved.
So assume that there exists some $b\in B\setminus\{0\}$ such that $b^2=0$.

By assumption, $B\subset \tilde B:=H^0(X,\Cal O_X)^G$, and
$\tilde B$ is a Noetherian $B$-module by what we have already proved.
Since $B$ is a $B$-submodule of $\tilde B$, we have that $B$ is a
Noetherian ring.
So we only need to prove that $B\red$ is of finite type.

Set $X_1$ to be the $G$-stable closed subscheme of $X$ defined by $b\Cal O_X$,
and let $Y_1$ be the scheme theoretic image of $X_1\hookrightarrow X
\rightarrow Y$.
Then $X_1\rightarrow Y_1$ is a universal 
strict orbit space, and hence $Y_1$ is of
finite type by induction assumption.
Hence $Y\red=(Y_1)\red$ is also of finite type, as desired.
\end{proof}

\begin{proof}[Proof of Theorem~\ref{geometric-quotient.thm}]
We only need to check the assumption of Lemma~\ref{reduction-lemma.thm}.
So we may assume that $Y$ is reduced, and it suffices to prove that 
$Y$ is of finite type.
We may assume that $S=\Spec R$ and $Y=\Spec B$ are affine.

Note that for any Noetherian flat $R$-algebra $R'$, the base change
$\varphi':X'\rightarrow Y'$ is again a universal strict orbit space.
%, since we
%know that $\varphi$ is universally open.
By Onoda's theorem \cite[(2.11) and (2.20)]{Onoda}, we may assume that
$R$ is local.
Since $\hat R$ is excellent, $\hat R\otimes_R Y$ is of finite type
again by Lemma~\ref{reduction-lemma.thm}.
By the descent argument \cite[(2.7.1)]{EGA-IV}, $Y$ is of finite type.
\end{proof}

\section{Proper morphisms}

Let $R$ be a Noetherian ring, and 
$\varphi:X\rightarrow Y$ a surjective $R$-morphism of $R$-schemes with
$X$ of finite type.

\begin{lmm}\label{universally-catenary.thm}
Let $R$ be universally catenary, $\varphi:X\rightarrow Y$ proper,
and $Y$ Noetherian.
Then $Y$ is universally catenary \(i.e., 
all local rings of $Y$ are universally catenary\).
\end{lmm}

\begin{proof}
Replacing $Y$ by an integral scheme of finite type over $Y$, 
it suffices to show that
under the assumption of the lemma, if $Y=\Spec B$ is affine and
integral, $Q,P\in\Spec B$ with $Q\subset P$, $\height(P/Q)=1$, then
$\height P=\height Q+1$.
Replacing $X$ by an irreducible component with the reduced structure
that is surjectively mapped onto $Y$, we may assume that $X$ is integral.

By Proposition~\ref{prop.thm},
there exists some finitely generated $R$-subalgebra $A\subset B$
such that $\eta:Y=\Spec B\rightarrow \Spec A$ is birational and 
geometrically injective.
Since $\varphi$ is surjective, the dimension formula holds
by \cite[(1.11)]{Onoda}, and we have $\height P=\height (P\cap A)$
and $\height Q=\height(Q\cap A)$.

As $\varphi$ is a closed morphism, there is a sequence 
$x_0,x_P,x_Q$ of points of $X$ such that 
$x_0$ is a closed point of the generic fiber,
$x_Q$ is a specialization of $x_0$ and $f(x_Q)=Q$, and
$x_P$ is a specialization of $x_Q$ and $f(x_P)=P$.
Then by the dimension formula, 
we have
\begin{eqnarray*}
\dim \Cal O_{X,x_P} & = & \height P +\tdeg_{R(Y)}R(X)-
                               \tdeg_{\kappa(P)}\kappa(x_P)\\
\dim \Cal O_{X,x_Q} & = & \height Q +\tdeg_{R(Y)}R(X)-
                               \tdeg_{\kappa(Q)}\kappa(x_Q)\\
\dim \Cal O_{\overline{x_Q},x_P} & = & \height(P/Q)+
                                 \tdeg_{\kappa(Q)}\kappa(x_Q)
                               -\tdeg_{\kappa(P)}\kappa(x_P).
\end{eqnarray*}
Since $\Cal O_{X,x_P}$ is catenary, we have that
\[
\dim \Cal O_{X,x_P}=\dim \Cal O_{X,x_Q}+\dim \Cal O_{\overline{x_Q},x_P}.
\]
Hence $\height P=\height Q+\height(P/Q)=\height Q+1$, as desired.
\end{proof}

\begin{thm}\label{proper.thm}
Let $R$ be an excellent ring, and $f:X\rightarrow Y$ a surjective proper
morphism of $R$-schemes.
If $X$ is of finite type over $R$ and $Y$ is a Noetherian scheme,
then $Y$ is of finite type over $R$.
\end{thm}

\begin{proof} We may assume that $Y=\Spec B$ is affine and integral.
We may assume that $X$ is integral.
As in the proof of Lemma~\ref{universally-catenary.thm}, 
we take a finitely generated $R$-subalgebra $A$ of $B$ such that
$\Spec B\rightarrow \Spec A$ is geometrically injective.
The dimension formula holds between $A$ and $B$, and $B$ is universally
catenary.
By \cite[(4.9)]{Onoda}, $B$ is of finite type.
\end{proof}

\end{document}